\documentclass[12pt]{article}

\usepackage{amsmath, amssymb, latexsym, amscd, amsthm} 

\newtheorem{theorem}{Theorem}[section]
\newtheorem{corollary}[theorem]{Corollary}
\newtheorem{lemma}[theorem]{Lemma}

\newtheorem{statement}[theorem]{Statement}

\theoremstyle{definition}

\newtheorem{remark}[theorem]{Remark}

\def\diam{\operatorname{diam}}
\def\id{\operatorname{id}}
\def\Int{\operatorname{Int}}
\def\Emb{\operatorname{Emb}}

\begin{document}

\title{All projections of a typical Cantor set are Cantor sets}   
\author{Olga Frolkina\footnote{Supported by Russian Foundation of Basic Research (Grant
No.~19--01--00169).}\\
Chair of General Topology and Geometry,\\
Faculty of Mechanics and Mathematics,\\
M.V.~Lomonosov Moscow State University,\\
Leninskie Gory 1, GSP-1,\\
Moscow 119991, Russia\\
E-mail: olga-frolkina@yandex.ru
}

\maketitle

\begin{abstract}
In 1994, John Cobb asked:
given $N>m>k>0$, does there exist a Cantor set in $\mathbb R^N$ such that
each of its projections
into $m$-planes is exactly $k$-dimensional?
Such sets
were described
for $(N,m,k)=(2,1,1)$ by L.~Antoine (1924) and for
$(N,m,m)$ by K.~Borsuk (1947).
Examples 
were
constructed for the cases
$(3,2,1)$ by J.~Cobb (1994),
for $(N,m,m-1)$ and in a different way for $(N,N-1,N-2)$ by O.~Frolkina (2010, 2019),
for $(N,N-1,k)$ by S.~Barov,
J.J.~Dijkstra and M.~van der Meer (2012).
We show that such sets are exceptional
in the following sense.
Let $\mathcal C(\mathbb R^N)$ be a 
set of all Cantor subsets of $\mathbb R^N$
endowed with the Hausdorff metric.
It is known that $\mathcal C(\mathbb R^N)$
is a Baire space.
We prove that there is
a dense $G_\delta $ subset $\mathcal P \subset
\mathcal C(\mathbb R^N)$ such that
for each $X\in \mathcal P$ and
each non-zero linear subspace $L
\subset \mathbb R^N$, the orthogonal projection
of $X$ into $L$ is a Cantor set.
This gives a partial answer 
to another question of J.~Cobb stated in the same paper (1994).


\end{abstract}

\section{Introduction and
Statements of Main Results}

Any topological space homeomorphic to 
the standard middle-thirds Cantor set 
$\mathcal C\subset I=[0,1]$ 
is called a Cantor set.
In 1884,  G.~Cantor
described a 
continuous surjection of $\mathcal C$ onto $I$;
it takes a point
$\frac{x_1}{3} + \frac{x_2}{3^2} + \ldots 
\in\mathcal C$ 
(where each
$x_i\in \{0;2\}$)
into the point
$\frac{x_1}{2}+\frac{x_2}{2^2} + \ldots \in I$ \cite[p.~255--256]{Cantor}.
In \cite{Zoretti},
L.~Zoretti
recalls 
as a known fact
(``ce resultat bien connu etant acquis...'')
that
the graph 
of any such surjection
is a Cantor set in plane 
whose projection to $y$-axis coincides with $I$. 
Taking an appropriate union of such sets,
he constructed
a Cantor set in plane
such that its projections
to a countable set of lines 
contain segments
\cite[p.~763]{Zoretti}.
Zoretti made an attempt
to describe a Cantor 
set in plane all of whose projections 
contain segments \cite[p.~763]{Zoretti} 
(see also \cite{Zoretti-ZB}); 
unfortunately it is not clear from the text why this set does not contain arcs.
(The purpose of Zoretti's note was to
disprove a statement made by F.~Riesz who wrote in \cite[p. 651]{Riesz}
that a projection of a closed
hereditarily disconnected subset of plane
to a straight line is again closed and
hereditarily disconnected.
F.~Riesz attributed this proposition to R.~Baire; this is incorrect because the work of Baire contained a different statement.)

L.~Antoine
described 
a Cantor set in $\mathbb R^2$
all of
whose projections coincide
with those of a regular hexagon
\cite[p.~272, {\bf 9}; p.~273, fig.~2]{Antoine-FM}.
Other examples 
of Cantor sets in plane all of whose projections
are segments
can be found in
\cite[p.~124, Example]{Cobb-projections},
\cite[Prop.~1]{Dijkstra-van-Mill}.

K.~Borsuk \cite{Borsuk} constructed a 
Cantor set in $\mathbb R^N$,  $N\geqslant 2$, such that its projection to every hyperplane contains an $(N-1)$-dimensional ball
(see \cite[Prop.~3.1]{DGW} for an alternative proof).
As a corollary, 
for any $m\leqslant N-1$
the projection of this set to each
$m$-plane
has dimension $m$.
Extensions of this result can be found in
\cite{Meyerson},
\cite{Glaser},
\cite{Kuzminykh}.
In contrast to the case of plane,
for $N\geqslant 3$ in $\mathbb R^N$
there is no Cantor set 
all of whose projections to $(N-1)$-planes
are convex bodies
\cite[Thm.~3]{Cobb-projections};
generalizations are obtained in
\cite[Thm.~4.7]{DGW},
\cite[Thm.~1]{BCD}.

J.~Cobb~\cite[Thm.~1]{Cobb-projections}  gave an example of a Cantor set in $\mathbb R^3$ whose projection to
every $2$-plane is $1$-dimensional.
He asked \cite[p.~126]{Cobb-projections}:
given $N>m>k>0$, does there exist a Cantor set in $\mathbb R^N$ such that
each of its projections
into $m$-planes is exactly $k$-dimensional?
(As in \cite{BDM}, call it the
case
$(N,m,k)$ of Cobb's problem.)
For the cases $(N,m,m-1)$ and $(N,N-1,k)$,
examples are given in
\cite[Thm.~1]{Frolkina-proj}
and 
\cite[Thm.~1]{BDM}, respectively.
I recently remarked that
a well-known Antoine's Necklace can serve as 
an example of a 
Cantor set in $\mathbb R^3$
all of whose projections are one-dimensional and 
connected  \cite{Frolkina-Arh80}; further,
each Cantor set in $\mathbb R^N$ can be moved 
by a small 
ambient isotopy
to an $(N,N-1,N-2)$-set
\cite{Frolkina-Gr}.

For a set $X\subset \mathbb R^N$,
Cobb suggested a notion of 
\emph{general position with respect to all projections} \cite[p.~127]{Cobb-projections}.
In some sense,
the behaviour of such sets 
under projections
is
as close as possible to that of embeddings;
all projections of such a Cantor set
are again Cantor sets.
In $\mathbb R^N$, $N\geqslant 2$,
Cobb
constructed 
a Cantor set in 
general position with respect to all projections
\cite[Thm.~5]{Cobb-projections}.
Cobb 
asked \cite[p.~128]{Cobb-projections}:
``Cantor sets that raise dimension under all projections
and those in general position with respect to all projections are 
both dense in the Cantor sets in $\mathbb R^m$ ---
which (if either) is more common, in the sense of
category or dimension or anything?''
Our main result provides a partial answer:
all projections
of a typical Cantor set in 
$\mathbb R^N$, $N\geqslant 1$,
are Cantor sets. 

\begin{theorem}\label{main}
For $N\geqslant 1$,
in $
\mathcal C(\mathbb R^N)$
there is
a dense $G_\delta $ subset $\mathcal P $ 
such that
for each $X\in \mathcal P$ and
each non-zero linear subspace $L
\subset \mathbb R^N$, the orthogonal projection
of $X$ into $L$ is a Cantor set.
\end{theorem}

Here, the set 
$\mathcal C(\mathbb R^N)$
of all Cantor sets in $\mathbb R^N$
is endowed with the Hausdorff metric;
it is 
a dense $G_\delta $ subset
of the space 
 $\mathcal K (\mathbb R^N)$
of all non-empty
compact subsets
of $\mathbb R^N$
(in \cite[p.~1089]{Gartside},
this fact is called ``well known'';
for a proof, see
\cite[Lemma 1]{Kuzminykh-typ} or \cite[Lemma 2.1, Remark 2.2]{Gartside}).
Thus $\mathcal C(\mathbb R^N)$ 
is a Baire space.

By the trivial projection we mean the projection into
the subspace $\{0\}$.

\begin{corollary}\label{space-of-compacta}
For $N\geqslant 1$,
 there is a dense $G_\delta $ subset $\mathcal A \subset  \mathcal K(\mathbb R^N)$ 
such that each $X\in \mathcal A$
is a Cantor set 
all of whose non-trivial projections are Cantor sets.
\end{corollary}

Denote by $\Emb (\mathcal C,\mathbb R^N)$
the space of all embeddings of the 
Cantor set $\mathcal C$
into $\mathbb R^N$ endowed  
with the distance
$\rho (f,g) = \sup \{ d(f(x),g(x)) \mid x\in\mathcal C\} $;
it is a Baire space
\cite[Thm. V~2]{HW}.
The evaluation map
$\mathcal E: \Emb (\mathcal C,\mathbb R^N)\to
\mathcal C(\mathbb R^N)$,
$f\mapsto f(\mathcal C)$
is an open continuous surjection
(Lemma \ref{openness});
hence Theorem \ref{main}
implies

\begin{corollary}\label{space-of-embeddings}
For $N\geqslant 1$,
in $ 
\Emb (\mathcal C,\mathbb R^N)$
there is a dense $G_\delta $ subset
$\mathcal B $ 
such that
for each $f\in\mathcal B$,
any non-trivial
projection
of $f(\mathcal C)$
is a  Cantor set.
\end{corollary}

Cantor sets are exactly
non-empty metric zero-dimensional
perfect compacta
\cite[7.4]{Kechris}.
Thus
Theorem \ref{main} follows 
from three 
statements
(here $p_L X$ denotes the orthogonal projection of $X$
into $L$):

\begin{statement}\label{one-point}
For $N\geqslant 1$,
let
$\mathcal D$ be
the set of all 
$X\in \mathcal C(\mathbb R^N)$
such that 
$|p_L X | = 1$
for some non-zero linear subspace
$L \subset \mathbb R^N$.
Then $\mathcal D$
is closed and nowhere dense in 
$\mathcal C(\mathbb R^N)$.
\end{statement}

\begin{statement}\label{isolated-point}
For $N\geqslant 1$, let $\mathcal I$
be the set of all 
$X\in \mathcal C(\mathbb R^N)$
such that 
for some non-zero linear subspace
$L \subset \mathbb R^N$ holds:
$|p_L X| > 1$ and $p_L X$ has an isolated point.
Then $\mathcal I$
is a meager $F_\sigma $ subset 
of
$\mathcal C(\mathbb R^N)$.
\end{statement}

\begin{statement}\label{zero-dim}
For $N\geqslant 1$,
there is
a dense $G_\delta $ subset $\mathcal Z \subset
\mathcal C(\mathbb R^N)$ such that
$\dim p_L X  =0 $
for each $X\in \mathcal Z$ and
each linear subspace $L
\subset \mathbb R^N$.
\end{statement}

To obtain Theorem \ref{main}, put
$\mathcal P =  \mathcal Z - \mathcal D - \mathcal I$.

\begin{remark}
We restrict ourselves by orthogonal projections only for simplicity of exposition and
notation.
In fact, suppose that
$\mathbb R^N = L_1 \oplus L_2$
and take any set $X\subset \mathbb R^N$.
The projection of $\mathbb R^N$
into $L_1$  parallel to $L_2$
maps $X$ to a set
affinely homeomorphic to the image
of $X$ under the orthogonal projection
of $\mathbb R^N$ into $L_2 ^\bot $.
\end{remark}

\section{Preliminaries}\label{Prelim-results}

Let $X$ be a non-empty topological space.
A set $A\subset X$ is called 
\emph{meager in $X$}, or \emph{of first category in $X$},
if it is the union of countably many
nowhere dense subsets of $X$.
A set $A$ that is not of first category in $X$ is said to be
\emph{non-meager in $X$} or
\emph{of second category in $X$}.
A set $A\subset X$ is called \emph{comeager in $X$}
or 
\emph{residual in $X$}
if it is the complement of a meager set;
equivalently, if
it contains the intersection
of a countable family of dense open sets.
If the space $X$ is of first category (in itself),
then each of its subsets is of first category
in $X$;
hence the two last definitions make sense only
if $X$ is of second category (in itself).

A \emph{Baire space}
is a topological space $X$ with the property:
every non-empty open set in $X$ is of second category in $X$;
equivalently, every comeager set is dense in $X$;
equivalently, 
for each countable family $\{ U_n , n\in\mathbb N \}$ 
of open dense subsets of $X$, their intersection
$\bigcap\limits_{n=1}^\infty U_n$
is dense in $X$.
Any Baire space is of
second category in itself.
\emph{Baire's Theorem} states that
\emph{a completely metrizable space is Baire}
(see e.g. \cite[8.4]{Kechris} for a proof).
A $G_\delta $ subset of a completely metrizable
space is itself completely metrizable
\cite[3.11]{Kechris}.

For a non-empty Baire space $X$
and its subset $P$,
we say that
\emph{a generic (or typical) element of $X$ is in $P$}
or 
\emph{most elements of $X$ are in $P$}
if $P$ is comeager in $X$;
this makes sense since in such a space,
it is impossible that
``most elements are in $P$''
and ``most elements are in $X-P$''
simultaneously.
Any dense $G_\delta $ subset of a Baire space $X$ is comeager in $X$.

\subsection*{Notation}

For a linear subspace $L$ of $\mathbb R^N$,
denote by $p_L$ the orthogonal projection
of $\mathbb R^N$ into $L$.
A projection $p_L$ is called non-trivial if 
$L\neq \{0\} $.

The Euclidean distance
between points $x,y\in \mathbb R^N$ is denoted by $d(x,y)$;
in general, for any
non-empty subsets $X,Y \subset \mathbb R^N$ we define
$d(X,Y) = \inf \{ d(x,y) \mid x\in X,\  y\in Y\} $.
For a non-empty set $X\subset \mathbb R^N$,
let $O(X, \varepsilon )=
\{ y\in\mathbb R^N  \mid d(y,X) < \varepsilon \}$.
In particular, $O(a,\varepsilon )$
is the open ball with center $a$
and radius $\varepsilon $.
By $B(a, \varepsilon )$ we denote 
the closed ball with center $a$
and radius $\varepsilon $.
(The dimension $N$ of the ambient space
is not reflected in the notation
$O(a,\varepsilon )$ and $B(a, \varepsilon )$;
it is usually clear from the context.)

The Hausdorff distance between
non-empty compacta $X,Y\subset \mathbb R^N$
is $d_H (X,Y) = 
\inf \{ \varepsilon > 0 \mid
X\subset O(Y,\varepsilon )
\ \text{and} \ 
Y\subset O(X,\varepsilon )\}$.
(Equivalently, 
$d_H (X,Y) < \varepsilon $
iff
$X\subset O(Y, \varepsilon )$
and 
$Y\subset O(X, \varepsilon )$.)

By $G(\ell , \mathbb R^N)$
we denote the Grassman space
consisting
of all $\ell $-dimensional linear subspaces of 
$\mathbb R^N$;
the distance 
between $L_1$ and $L_2$
is defined
as the Hausdorff distance 
$d_H( L_1\cap S, L_2\cap S)$, where $S\subset \mathbb R^N$ is the unit sphere
\cite{Gohberg-Markus}.
The space $G(\ell , \mathbb R^N)$ is compact
\cite[Cor. 6.2, Cor. 6.7]{Berkson}.

The space 
$\mathcal K (\mathbb R^N)$
of all non-empty
compact subsets
of $\mathbb R^N$
is endowed with the Hausdorff metric
(equivalently, one may use the Vietoris topology
\cite[4.F, 4.21]{Kechris});
this space is complete
\cite[4.25]{Kechris}.

Any topological space
homeomorphic to the usual middle-thirds Cantor set 
is called a Cantor set.
(These are exactly non-empty metric zero-dimensional perfect compacta
\cite[7.4]{Kechris}.)

The set
$\mathcal C(\mathbb R^N)$
of all Cantor sets in $\mathbb R^N$
is 
a dense $G_\delta $ subset of  $\mathcal K (\mathbb R^N)$
(in \cite[p.~1089]{Gartside},
this fact is called ``well known'';
the proof can be found in
\cite[Lemma 1]{Kuzminykh-typ} or \cite[Lemma 2.1, Remark 2.2]{Gartside}).
Hence $\mathcal C(\mathbb R^N)$ is a Baire space.
Although the Hausdorff metric
on $\mathcal C(\mathbb R^N)$
is not complete,
we may --- and will --- use it to prove
purely topological statements (concerning e.g. 
openness, nowhere density, etc.).

By $\Emb (\mathcal C,\mathbb R^N)$ we denote
the space of all embeddings of the 
Cantor set $\mathcal C$
into $\mathbb R^N$,
with the distance
$\rho (f,g) = \sup \limits_{x\in\mathcal C} d(f(x),g(x))$;
it is a Baire space.

For a topological manifold-with-boundary
$M$,
denote by 
$\Int M$
and $\partial M$
the interior and the boundary of $M$, correspondingly.

A system of points
$x_0, x_1, \ldots  , x_m$ of $\mathbb R^N$
is said to be \emph{in general position},
if each of its subsystems
of $k+1$ points $\xi _0,\xi_1,\ldots , \xi _k $ 
($k\leqslant N$) is affinely independent,
see \cite{Pontryagin-ct}.
Any finite system of points of $\mathbb R^N$
can be brought into general position
by an arbitrarily small displacement
\cite[Thm.~2]{Pontryagin-ct}.
For any finite system of points 
$x_0, x_1, \ldots  , x_m$ 
of $\mathbb R^N$ in general position,
there exists an $\varepsilon >0$
such that
$d(x_i,y_i)<\varepsilon $, $i=0,\ldots , m$
implies that the system
$y_0, y_1, \ldots  , y_m$ 
is in general position \cite[Thm.~1]{Pontryagin-ct}.

Finally, $|A|$ denotes the cardinality of a set $A$.

\section{Proofs of Statements \ref{one-point}
and \ref{isolated-point}} 

We start with two useful lemmas.

\begin{lemma}\label{prelim}
Let $N\geqslant 1$.
1) For each 
linear subspace $L\subset \mathbb R^N$
and each pair of non-empty compact sets $X, Y \subset \mathbb R^N$
we have
$d_H(p_L X , p_L Y) \leqslant d_H (X,Y)$.
\\
2) Let $\lim\limits _{n\to \infty} L_n = L$ in the space $G(\ell , \mathbb R^N)$.
Then
$\lim\limits_{n\to \infty } d_H (p_{L_n} X , p_L X)= 0 $ 
for each non-empty compactum $X \subset \mathbb R^N$.
\\
3) For each
$\ell = 0,\ldots  ,N$
the map
$$
\mathcal K (\mathbb R^N) \times G(\ell , \mathbb R^N)
\to 
\mathcal K (\mathbb R^N) , 
\quad
(X,L) \mapsto p_L X
$$
is continuous.
\end{lemma}

Proof.
1) 
For each $x,y\in\mathbb R^N$
holds
$d(p_L x, p_L y)  \leqslant d(x,y)$;
this implies the proposition.

2) The statement can be easily verified when 
$X$ is
a closed 
$N$-dimensional 
ball $B(a,r)$.
In general case, for each 
$\varepsilon >0$
take a finite family 
$B_1,\ldots , B_t$
of closed balls
 such that
$X\subset \cup_{i=1}^t B_i \subset 
O\left(  X, {\frac{\varepsilon }{2}}\right)$.
There exists an integer $n_0$ such that 
$d_H (p_{L_n} B_i , p_L B_i ) <
\frac{\varepsilon }{2}$
for each
$n\geqslant n_0$ and each
$i=1,\ldots ,t$.
Hence
$p_{L_n } B_i \subset 
O\left( p_L B_i, {\frac{\varepsilon }{2}} \right)$
and
$p_{L} B_i \subset 
O\left( p_{L_n} B_i, {\frac{\varepsilon }{2}}\right)$;
this implies
$$
p_{L_n} X \subset
p_{L_n} \left( \bigcup\limits_{i=1}^t B_i\right)
=
 \bigcup\limits_{i=1}^t p_{L_n} B_i
 \subset
 \bigcup\limits_{i=1}^t O\left( p_{L} B_i,
 {\frac{\varepsilon }{2}}\right) =
O\left( \bigcup\limits_{i=1}^t p_{L} B_i,
{\frac{\varepsilon }{2}} \right) =
$$
$$
=
O \left( p_{L} \left( \bigcup\limits_{i=1}^t B_i\right), {\frac{\varepsilon }{2}}\right) \subset
O \left(p_{L} \left( O \left( X ,{\frac{\varepsilon }{2}}\right)\right) , {\frac{\varepsilon }{2}}\right) \subset
O(p_{L} X , {\varepsilon } ) .
$$
Similarly,
$
p_L X\subset O  (p_{L_n } X, \varepsilon)$.
We get
$d_H (p_{L_n } X , p_L X) < \varepsilon  $.

3) Suppose that $\lim\limits_{n\to \infty}
(X_n, L_n) = (X,L)$ in 
$\mathcal K (\mathbb R^N) \times G(\ell , \mathbb R^N)$.
We have
$$d_H(p_{L_n} X_n , p_L X) \leqslant
d_H (p_{L_n} X_n , p_{L_n} X ) +
d_H (p_{L_n} X , p_{L} X ) 
\leqslant 
$$
$$
\leqslant 
d_H ( X_n ,  X ) +
d_H (p_{L_n} X , p_{L} X ) ;
$$
this tends to zero as $n\to \infty  $.

\begin{lemma}\label{finite-appr}
Let $N\geqslant 1$.
For
any $\varepsilon >0$
and any non-empty compactum
$K\subset \mathbb R^N$
there exists a finite set $A\subset \mathbb R^N$
in general position
such that
$d_H(K,A) < \varepsilon $
and
$|A|\geqslant N+1$.
\end{lemma}

Proof.
 Let $F\subset K$
be 
 a finite $\frac{\varepsilon }{2}$-net
 of $K$;
we may assume that $|F|\geqslant N+1$. 
Moving each point of $F$
at a distance
less than $\frac{\varepsilon }{2}$
from its original position,
    construct
    a set
 $A\subset \mathbb R^N$
such that $|F|=|A|$,
 $A$ is in general position, 
and $d_H(A,F)<\frac{\varepsilon }{2}$.
We have
$d_H(K,A)\leqslant 
d_H (K,F) + d_H(F,A) < \varepsilon $;
thus $A$ is the desired set.

Proof of Statement \ref{one-point}.
\emph{Closedness.}
Take
$X_1,X_2, \ldots , X \in \mathcal C(\mathbb R^N)$
such that
$d_H(X_n,X) \to 0 $ as $n\to \infty $.
For each $n\in\mathbb N$
let $L_n \subset \mathbb R^N$ be an
$\ell $-dimensional linear subspace
such that
$|p_{L_n} X_n| = 1$.
The Grassman space
$G(\ell , \mathbb R^N)$
is compact;
we may  therefore assume
(passing to a subsequence if necessary)
that 
there exists an $L\in G(\ell , \mathbb R^N)$ such that
$L_n \to L$ as $n\to \infty $.
By 3) of Lemma \ref{prelim}
$\lim\limits_{n\to \infty } 
d_H (p_{L_n} X_n , p_L X) = 0$.
Recall that $|p_{L_n} X_n | =1$;
we easily get $|p_L X |=1$.

\emph{Nowhere density.}
Take any $X\in \mathcal C(\mathbb R^N)$
and any $\varepsilon >0$.
By Lemma~\ref{finite-appr}
there exists a finite set $A = \{ a_1,\ldots , a_t \} \subset \mathbb R^N$ 
in general position
such that
$d_H(X,A)<\frac{\varepsilon }{2}$
and
$t = |A| \geqslant N+1$.

Take a number $r>0$
such that:
$r<\frac12 \min \{ d(x,y ) \mid  x,y\in A; x\neq y \}$; $r<\frac{\varepsilon }{2}$;
and for any choice of points
$b_1\in O (a_1 ,r), \ldots , b_t\in O (a_t ,r)$
the set $\{ b_1,\ldots , b_t\}$
is in general position.

For each $i=1,\ldots ,t $
take a Cantor set
$K_i \subset O (a_i ,r)$.
Let $K=K_1\cup \ldots \cup K_t$;
this is a Cantor set.
We have
$d_H(K,A) < r < \frac{\varepsilon }{2}$,
hence
$d_H(X,K)\leqslant d_H(X,A)+d_H(A,K) < \varepsilon $.

Choose a positive number $\delta $ such that
$\delta < \min \limits_{i=1,\ldots , t}
d\left(K_i, \partial ( O (a_i ,r))\right)$.
Let $Y\subset \mathbb R^N$
be any  Cantor set 
with $d_H(Y,K)<\delta $.
Let us show that
$Y\notin \mathcal D$.
Suppose that $|p_L Y | = 1$
for some non-zero linear subspace $L\subset \mathbb R^N$.
Hence $Y$ lies in an affine subspace
$\vec v + L^\bot $ of $\mathbb R^N$.
For each $i=1,\ldots , t$
take a point $b_i \in O (a_i ,r) \cap Y$ 
(by the construction, this intersection is non-empty). 
Recall that the set $\{ b_1,\ldots , b_t\}$
is in general position and $t\geqslant N+1 \geqslant 
\dim (\vec v + L^\bot ) +2 $.
This contradiction proves Statement \ref{one-point}.

Proof of Statement \ref{isolated-point}.
For each $k\in \mathbb N$
define $\mathcal I_k$
as the set of all
$X\in \mathcal C(\mathbb R^N)$
such that
for some non-zero linear subspace $L\subset \mathbb R^N$
we have: $|p_L X|>1$ and
$d(p_L x , p_L X - p_L x) \geqslant \frac{1}{k}$
for some point $x\in X$.
Clearly $\mathcal I = \bigcup\limits_{k\geqslant 1} \mathcal I_k$.
We will prove the statement
by showing that
each $\mathcal I_k$ is closed and nowhere dense in $\mathcal C(\mathbb R^N)$.

\emph{Closedness.}
Suppose that
$X_1,X_2,\ldots \in \mathcal I_k$,
$X\in\mathcal C(\mathbb R^N)$,
and $ d_H (X_n,X)\to 0 $ as $n\to \infty $.
For each $n\in\mathbb N$,
there exist a subspace $L_n \subset \mathbb R^N$ and a point $x_n\in X_n$
such that
$|p_{L_n} X_n| > 1$
and
$d(p_{L_n} x_n  , p_{L_n }X_n - p_{L_n} x_n ) 
\geqslant \frac{1}{k} $.
We may assume that $\dim L_n = \ell $ is the same for
each $n$; clearly $\ell \neq 0,N$.

Since the Grassman space $G(\ell , \mathbb R^N)$
is compact,
we may assume that there exists
$\lim\limits_{n\to \infty  } L_n =  L$.
We may also assume that the sequence $\{ x_n \}_{n\in\mathbb N}$
converges to some 
$x\in X$ as $n\to\infty $. 
Let us show that
$|p_L X| >1$
and
$d(p_L x, p_L X - p_L x) \geqslant \frac{1}{k}$.

By Lemma \ref{prelim},
the equalities 
$\lim\limits_{n\to \infty  } d_H (X_n,X) = 0$ and $\lim\limits_{n\to \infty  } L_n = L$ 
together
imply 
$\lim\limits_{n\to \infty  } 
d_H (p_{L_n} X_n , p_L X ) =0$.
Recall that $\diam p_{L_n} X_n  \geqslant \frac{1}{k}$ for each $n$;
hence $|p_L X | > 1$.

Assume that
there exists a point $y\in X$
such that
$0<d(p_L x, p_L y) < \frac{1}{k}$.
Take a positive number $\varepsilon $
such that
$$\varepsilon < \frac12 d(p_L x, p_L y)
\quad
\text{and}
\quad
\varepsilon < \frac12 \left( \frac{1}{k} - d(p_L x, p_L y)
\right).$$
Since $x_n \to x$ and $L_n\to L$ as $n\to\infty $,
we have $\lim\limits_{n\to\infty } d(p_{L_n} x_n , p_L x)  =0$ by Lemma \ref{prelim};
and there exists an integer $n_1$ such that 
$d(p_{L_n} x_n , p_L x) < \varepsilon $
for
each $n\geqslant n_1$.
Passing to a subsequence if necessary,
we may assume that for each $n\in \mathbb N$
a point $y_n\in X_n$ is chosen 
so that
$\lim\limits_{n\to\infty } y_n = y$;
in particular, there is an integer $n_2$
such that 
$d(p_{L_n} y_n , p_L y) < \varepsilon $
for each $n\geqslant n_2$
(see Lemma \ref{prelim}).
Now take any $n\geqslant \max \{ n_1,n_2\}$;
we have
$$
d(p_{L_n} x_n , p_{L_n} y_n )
\leqslant
d(p_{L_n} x_n , p_{L} x )
+
d(p_{L} x , p_{L} y )
+
d(p_{L} y  , p_{L_n} y_n ) <
$$
$$
<
\varepsilon +
d(p_{L} x , p_{L} y )
+ \varepsilon <\frac{1}{k} ,
$$
and
$$
d(p_{L_n} x_n , p_{L_n} y_n )
\geqslant
d(p_{L} x , p_{L} y )
-
d(p_{L_n} x_n , p_{L} x )
-
d(p_{L} y  , p_{L_n} y_n ) >
$$
$$
>
d(p_{L} x , p_{L} y )
- \varepsilon -
 \varepsilon >0.
 $$
This contradicts with the inequality 
 $d(p_{L_n}x_n, p_{L_n} X_n - p_{L_n} x_n) \geqslant \frac{1}{k}$.
Hence $d (p_L x, p_L X-p_L x)\geqslant \frac{1}{k} $;
that is, $X\in \mathcal I_k$.

\emph{Nowhere density.}
Fix an integer $k\geqslant 1$.
Take any $X\in \mathcal C(\mathbb R^N)$
and any $\varepsilon >0$.
We will construct a Cantor set 
$K\subset\mathbb R^N$ and a number $\delta >0$
such that
$d_H (X,K)<\varepsilon $
and 
$Y\notin \mathcal I_k$
for each Cantor set $Y\subset\mathbb R^N$
with 
$d_H(K,Y) < \delta $.

Take a finite set $A= \{ a_1,\ldots , a_t\}
\subset\mathbb R^N$
such that
$d_H(X,A) < \frac{\varepsilon }{2}$
(Lemma~\ref{finite-appr}).

Choose a positive number $r$ such that
$$r<\frac12 \min \{ d(x,y) 
 \mid x,y\in A ;\ x\neq y \};
 \quad r<\frac{\varepsilon }{2};
 \quad
 \text{and}
 \quad r<\frac{1}{2k}.$$

For each $i=1,\ldots , t$
take  a set
$A_i = \{ a_{i,1} , \ldots , a_{i,N+1} \}
\subset O(a_i , r)$
in general position.

Take a positive number $\rho $ such that
$$\rho < \frac12 \min \limits_{i=1,\ldots , t} \{ d( x , y) 
 \mid x,y\in A_i ;\
x\neq y\}
\quad
\text{and}
\quad
\rho <
\min \limits_{i=1,\ldots , t} \{
d(A_i, \partial ( O (a_i,r) ) ) \},
$$
and such that
for each $i=1,\ldots , t$
any set of points
$\{ b_{i,1} , \ldots , b_{i,N+1}\}$,
where
$b_{i,1} \in O(a_{i,1} ,\rho ) $,...,
$b_{i,N+1} \in O (a_{i,N+1}, \rho )$,
is in general position.
Note that $\rho < r <\frac{\varepsilon }{2}$.

For each $i=1,\ldots ,t$ and each
$j=1,\ldots , N+1$
take a Cantor set
$K_{i,j} \subset O(a_{i,j}, \rho )$.
The union $K=\bigcup \{
K_{i,j} \mid
i=1,\ldots ,t; \  j=1,\ldots , N+1\} $ 
is a Cantor set, and
$d_H (X,K) \leqslant d_H (X,A) + d_H (A,K)
<\frac{\varepsilon }{2} + \rho < \varepsilon $.

Take a positive number $\delta $ such that
$$
\delta < \min_{\substack{i=1,\ldots , t
\\
j=1,\ldots , N+1}} 
 d(K_{i,j} , \partial \left( O (a_{i,j} , \rho)\right) ) .
$$

Now let $Y$ be a Cantor set in $\mathbb R^N$
with $d_H(Y,K)<\delta $.
Let us show that $Y\notin \mathcal I_k$.
For each $i=1,\ldots , t$ and each
$j=1,\ldots , N+1$
the set
$Y_{i,j} =  O( K_{i,j} , \delta ) \cap Y$
is non-empty. 
Note that $Y_{i,j}\subset 
O(K_{i,j} , \delta ) \subset
O (a_{i,j} , \rho ) \subset O(a_i, r)$ for each $i$.
Let $L\subset \mathbb R^N$
be a linear subspace,
$\dim L \neq 0,N$. 
For each vector $\vec v\in \mathbb R^N$ and each $i=1,
\ldots , t$,
the affine subspace
$\vec v + L^\bot $
intersects at most $N$
balls of the family
$\{ O (a_{i,1}, \rho ) , \ldots , O(a_{i, N+1} , \rho )\}$.
In particular,
$|p_L (Y_{i,1} \cup 
\ldots \cup Y_{i,N+1})| \geqslant 2$;
hence $|p_L Y | \geqslant 2$.
For each $i=1,
\ldots , t$ we have
$$
\diam p_L (Y_{i,1} \cup \ldots \cup Y_{i,N+1})
\leqslant\diam
(Y_{i,1} \cup \ldots \cup Y_{i,N+1})
\leqslant
$$
$$
\leqslant
\diam
(O ( K_{i,1} , \delta )\cup \ldots \cup O( K_{i,N+1}, \delta ))
\leqslant \diam
O (a_i ,r) =2r <\frac{1}{k}.
$$
Thus
 $Y\notin \mathcal I_k$, and
the statement is proved.

\section{{Proof of Statement \ref{zero-dim}}}

We need some auxiliary results.

\begin{lemma}\label{lemma-1}
Let $\{ p_1,\ldots , p_q \}$ be 
 a set of pairwise distinct points in
$\mathbb R^D$, ${D\geqslant 1}$.
Suppose that
the union of 
$D$-dimensional
balls
$B( p_1 ,\delta )\cup \ldots \cup B( p_q , \delta )$ 
with radii $\delta $
is
connected and
$\diam (B( p_1 , \delta )\cup \ldots \cup B(p_q , \delta  ))
\geqslant 2\delta (N+1)$,
where $N$ is a positive integer.
Then there exists an integer $s\geqslant N+1$ and
an ordered set $(i_1, \ldots , i_s)$
of pairwise distinct integers such that:
$i_j \in \{ 1,\ldots , q\}$ for each $j=1,\ldots , s$; and
$d(p_{i_j} , p_{i_{j+1}}) \leqslant 2\delta $ for each
$j=1,
\ldots , s-1$.
\end{lemma}

Proof.
By assumption, there exist points
$x,y\in \bigcup\limits_{i=1}^q B( p_i , \delta )$
such that $d(x,y) \geqslant 2\delta (N+1)$.
Take any $i_1, \alpha \in \{ 1,\ldots , q\}$ such that
$x\in B(p_{i_1} , \delta  )$ and 
$y\in B(p_{\alpha } , \delta  )$;
note that $i_1$ and $\alpha $ can not coincide.
Moreover, we have
$$
d(p_{i_1}, p_{\alpha })
\geqslant
d(x,y) - d(x,p_{i_1}) - d(y,p_{\alpha })
\geqslant 2\delta (N+1) - \delta  -\delta = 2\delta N.
$$
Take the sequence of sets
$$\mathcal A_0 = \{ \alpha \};$$
$$\mathcal A_1 = \{ 
j\in \{ 1,\ldots , q\} \mid 
B( p_j , \delta ) \cap B( p_{\alpha } , \delta) \neq\emptyset 
\} ;$$
$$
\mathcal A_2 = \{
j\in \{ 1,\ldots , q\} \mid 
\exists m \in \mathcal A_1 : 
B(p_j ,\delta ) \cap B(p_{m} ,\delta )\neq\emptyset 
\} ;
$$
$$
\mathcal A_3 = \{
j\in \{ 1,\ldots , q\} \mid 
\exists m \in \mathcal A_2 : 
B( p_j , \delta )\cap B ( p_{m} , \delta  ) \neq\emptyset 
\} ;
$$
and so on.
We have $\mathcal A_0 \subset \mathcal A_1\subset \mathcal A_2\subset \ldots $.
Since $q$ is finite, for some $n$
we have $\mathcal A_n = \mathcal A_{n+1} = 
\mathcal A_{n+2} = \ldots $.
Connectedness of
$\bigcup\limits_{i=1}^q B(p_i , \delta  )$
implies that $\mathcal A_n = \{1,\ldots , q\}$.

The index $i_1$ is already defined; 
we extend it to an ordered set of
pairwise distinct
indices
$i_1,\ldots  ,i_s\in \{1,\ldots , q\}$
such that
$B(p_j , \delta  ) \cap B( p_{j+1} , \delta )\neq\emptyset $
for each $j=1,\ldots , s-1$.
This can be done as follows.
Take the minimal integer $m_1$
such that
$i_1 \in \mathcal A _{m_1}$.
By the definition of
$\mathcal A _{m_1}$,
there exists 
an $i_{2} \in \mathcal A _{m_{1}-1}$
such that
$B( p_{i_{1}} , \delta ) \cap B(p_{i_{2}} , \delta  )\neq\emptyset $. By the minimality of $m_1$ we have
$i_{2} \neq i_1$.
Let $m_{2}$ be the minimal possible
integer $m$ with the property
$i_{2} \in \mathcal A _{m}$.
Continue in similar way.
Find an index $i_{3} \in \mathcal A_{m_{2} - 1}$
such that
$B( p_{i_{2}} ,\delta )\cap B( p_{i_{3}} , \delta )\neq\emptyset $; the minimality property of $m_{2}$ and $m_1$
implies that $i_{3} \neq i_{2}$
and $i_{3}\neq i_1$. 
Denote by $m_{3}$ the 
 minimal possible integer $m$
such that 
$i_{3} \in \mathcal A _{m}$.
This process continues until
there will be no more possibility to lower the index;
in other words, there exists an integer $s$ such that
$i_s = \alpha $ and $i_j\neq\alpha $ for $j<s$.

The ordered set $(i_1,\ldots , i_s)$
has the desired properties.
Consider the broken line 
with subsequent vertices
$p_{i_1} , \ldots , p_{i_s}$. 
By construction, the length of each edge is
$d (p_{i_j} , p_{i_{j+1}}) \leqslant 2\delta $.
We have
$$
2\delta (s-1) \geqslant 
\sum \limits_{j=1}^{s-1} d (p_{i_j} , p_{i_{j+1}})
\geqslant d (p_{i_1} , p_{i_{s}})
\geqslant 2\delta N;
$$
therefore
$s\geqslant N+1$.
The lemma is proved.

{\bf Notation.}
Let $A
\subset \mathbb R^N$
be a finite set 
in general position, where ${N\geqslant 1}$.
Denote by
$\lambda (A)$
the infimum
of the sums
$\sum\limits_{j=1}^{N} d(p_L a_{i_j } , p_L a_{i_{j+1}})$
over all 
li\-ne\-ar subspaces $L\subset \mathbb R^N$,
$\dim L \neq 0, N$,
and all
ordered $(N+1)$-tuples
$(a_{i_1} , a_{i_2} , 
\ldots , a_{i_{N+1}})$ of pairwise distinct points
of $A$.

\begin{lemma}\label{lemma-3}
Let 
$A\subset \mathbb R^N$
be a finite set 
in general position,
$N\geqslant 1$,
and $|A|\geqslant N+1$.
Then $\lambda (A)>0$.
\end{lemma}

Proof.
It suffices to fix an arbitrary
ordered $(N+1)$-tuple 
$(a_{i_1}, a_{i_2} , \ldots , a_{i_{N+1}})$
of 
pairwise distinct
points of $A$ and an arbitrary
$\ell \neq 0,N$ and show that
the infimum of the sums
$\sum\limits_{j=1}^{N} d(p_L a_{{i_j }} , p_L a_{i_{j+1}})$
over all
$\ell $-dimensional
linear subspaces $L\subset \mathbb R^N$
is positive.
To prove this, consider
the map
$\Phi  _{\ell } : G(\ell ,\mathbb R^N) \to \mathbb R$
which takes
$L
\mapsto
\sum\limits_{j=1}^{N} d(p_L a_{j } , p_L a_{j+1})$;
it is continuous.
Recall that $A$
is in general position;
hence $\Phi _{\ell } (L) > 0$ for each $L \in G(\ell ,\mathbb R^N)$.
Since $G(\ell ,\mathbb R^N)$ is compact,
the lemma follows.

\begin{lemma}\label{lemma-2}
Let 
$A=\{ a_1,\ldots , a_t\} \subset \mathbb R^N$
be a set 
in general position, 
$N\geqslant 1$,
and $t=|A|\geqslant N+1$.
Let $\delta $ be a positive number
such that
$\delta < \frac{\lambda (A)}{2(|A| -1)} $.
Then for any linear subspace $L\subset \mathbb R^N$
the diameter of
each connected component
of
$p_L (B (a_1 , \delta ) \cup \ldots \cup B(a_t, \delta )) $
is less than $2\delta (N+1)$.
\end{lemma}

Proof.
Suppose that there exist a subspace $L\subset \mathbb R^N$
and a connected component $Z$ of 
$p_L (B(a_1 , \delta )\cup \ldots \cup B(a_t , \delta  )) $
such that
$\diam Z \geqslant 2\delta (N+1)$.
Each $p_L (B(a_i , \delta ))$
is a closed ball in $L$ 
of radius $\delta $ and center
$p_L a_i$;
hence we may assume
that
$Z = p_L (B( a_1 , \delta )\cup \ldots \cup B( a_q,
\delta )) $
for some $q\leqslant t$.
Denote $p_i := p_L a_i$ for brevity.
By Lemma \ref{lemma-1} there exists
an ordered set of pairwise distinct indices
$(i_1,\ldots , i_s)$ which take values
in $\{ 1,\ldots , q\}$
such that
$s\geqslant N+1$
and $d(p_{i_j} , p_{i_{j+1}}) \leqslant 2\delta $
for each $j=1,\ldots , s-1$.
We have
$$
\sum\limits_{j=1} ^{s-1} d(p_{i_j} , p_{i_{j+1}})
\leqslant 2\delta (s-1) \leqslant 2\delta (t-1)
< 
\lambda (A).
$$
On the other hand,
the inequality $s\geqslant N+1$ and the definition
of $\lambda (A)$ imply that
$$
\sum\limits_{j=1} ^{s-1} d(p_{i_j} , p_{i_{j+1}})
\geqslant \lambda (A).
$$
This contradiction proves the lemma.

Proof of Statement \ref{zero-dim}.
For each $k\in \mathbb N$,
let
$\mathcal Z _k$
be the collection of all
Cantor sets $X\subset \mathbb R^N$
such that
there exists
a finite family
$B_1,\ldots , B_s$
of pairwise disjoint 
closed $N$-dimensional balls
satisfying the following conditions:

(1) 
$X\subset \bigcup\limits_{i=1}^s \Int B_i$;

(2)
$X\cap \Int B_i \neq\emptyset $ for each $i=1,\ldots , s$;

(3)
$\diam B_i < \frac{1}{k}$ for each $i=1,\ldots , s$;

(4)
for any linear subspace
$L\subset\mathbb R^N$ 
the diameter of each connected component of
$p_L \left(\bigcup\limits_{i=1}^s B_i\right)$
is less than $\frac{1}{k}$.

We define the required subset $\mathcal Z$ as the intersection 
$\bigcap\limits_{k=1}^\infty
\mathcal Z_k$.
Then 
$\dim p_L X = 0$ for any $X\in \mathcal Z$ 
and
any linear
subspace $L\subset \mathbb R^N$.

Let us fix an arbitrary $k\in\mathbb N$
and prove that 
 $\mathcal Z_k$
is open and dense 
in 
$\mathcal C(\mathbb R^N)$.

\emph{Openness.}
This can be easily verified directly.
Alternatively,
note that for each 
family $B_1,\ldots , B_s$
the set of all $X\in \mathcal C(\mathbb R^N)$
which satisfy conditions (1) and (2)
is open in Vietoris topology.
Recall that
our
topology on $\mathcal C(\mathbb R^N)$
coincides with the
Vietoris topology.
(Conditions (3) and (4) simply put restrictions
on admissible balls $B_1,\ldots , B_s$.)

\emph{Everywhere density.}
Take any $K\in \mathcal C(\mathbb R^N)$
and any $\varepsilon >0$.
We will construct a Cantor set $X\subset \mathbb R^N$
such that
$X\in\mathcal Z_k$ and
$d_H(X,K)<
\varepsilon $.

Using Lemma \ref{finite-appr},
take 
a finite set $A=\{a_1,\ldots  , a_t\}$
in general position
such that
$t=|A|\geqslant N+1$
and
$d_H(K,A) < \frac{\varepsilon }{2}$.
Choose a positive number $\delta $ such that
$$\delta < \frac12 \min \{  
d(x , y)
\mid 
x,y\in A ;\ x\neq y
\} ;$$
$$
\delta < \frac{\lambda (A)}{2(|A|-1)} ;
\quad 
\delta < \frac{\varepsilon }{2} ; 
\quad
\delta < \frac{1}{2k(N+1)} .
$$
By Lemma \ref{lemma-2}
for each linear subspace $L\subset \mathbb R^N$
and each connected component $Z$
of 
$p_L (B(a_1 , \delta ) \cup \ldots \cup
B( a_t , \delta ))$
we have: $\diam Z < 2\delta (N+1) < \frac{1}{k} $.
For each $i=1,\ldots , t$
take a Cantor set $X_i \subset O(a_i , \delta )$.
The union $X=X_1\cup \ldots \cup X_t$ is the required Cantor set.
In fact,
let 
the balls $B(a_1 , \delta ), \ldots , 
B(a_t , \delta )$
play the role of $B_1,\ldots , B_s$ from the definition
of $\mathcal Z_k$;
we get
$X\in \mathcal Z_k$.
Finally,
our construction implies 
$d_H(X,A)<\delta <\frac{\varepsilon }{2} $
and hence
$d_H (X,K)\leqslant d_H(X,A)+d_H(A,K) < \frac{\varepsilon }{2} + \frac{\varepsilon }{2} = \varepsilon $.
The statement is proved.

\section{Appendix}

The proof of Theorem \ref{main}
is rather complicated.
Let us give a short argument for the 
easier case:
when the set of projections is at most countable.

\begin{statement}\label{easier}
Let
$L \subset \mathbb R^N$
be a linear subspace,
$N\geqslant 2$, and
$\dim L \neq 0,N$.

1) The set
$\mathcal P_1= \{
X\in \mathcal C(\mathbb R^N) 
\mid
p_L X \text{ is a Cantor set} \}$
is a dense $G_\delta $ subset of 
$ \mathcal C(\mathbb R^N) $.

2)
The set
$\mathcal P_2 = \{
X\in \mathcal K(\mathbb R^N) 
\mid
p_L X \text{ is a Cantor set} \}$
is a dense $G_\delta $ subset of 
$ \mathcal K (\mathbb R^N) $.

3)
The set
$\mathcal P_3 = \{
f \in \Emb ( \mathcal C , \mathbb R^N) 
\mid
p_L (f(\mathcal C )) 
\text{ is a Cantor set} \}$
is a dense $G_\delta $ subset of 
$ \Emb ( \mathcal C , \mathbb R^N)  $.
\end{statement}

This immediately implies

\begin{corollary}
For 
$N\geqslant 2$,
let 
$L _1, L _2,\ldots \subset \mathbb R^N$
be at most countable family of linear subspaces
such that
$\dim L _k \neq 0, N$ for each $k\in\mathbb N$.
The sets
$$\{
X\in \mathcal C(\mathbb R^N) 
\mid
p_{L_k} X \text{ is a Cantor set for each } k \},$$
$$\{
X\in \mathcal K (\mathbb R^N) 
\mid
p_{L_k} X \text{ is a Cantor set for each } k \}$$
and
$$\{
f \in \Emb ( \mathcal C , \mathbb R^N) 
\mid
p_{L_k} (f(\mathcal C )) 
\text{ is a Cantor set for each } k \}$$
are dense $G_\delta $ subsets
of
$\mathcal C(\mathbb R^N) $,
$\mathcal K (\mathbb R^N) $
and 
$\Emb ( \mathcal C , \mathbb R^N) $
correspondingly.
\end{corollary}

Let us first prove

\begin{lemma}\label{openness}
For each $N\geqslant 1$,
the map $\mathcal E :\Emb (\mathcal C,\mathbb R^N)
\to \mathcal C(\mathbb R^N)$,
$f\mapsto f(\mathcal C)$
is a 
continuous open surjection.
\end{lemma}

Proof.
Surjectivity and continuity of
$\mathcal E$ is clear.
Let us prove that $\mathcal E$ is an open map.
Fix any $f\in \Emb (\mathcal C,\mathbb R^N)$
and any $\varepsilon >0$.
It suffices to show that
$f(\mathcal C)$
is an interior point
of $\mathcal E ( U (f,\varepsilon ))$ in $\mathcal C(\mathbb R^N)$,
where
$U (f,\varepsilon )
= \{ g\in \Emb (\mathcal C,\mathbb R^N)\mid
\rho (g,f) < \varepsilon \}$.
We will find a $\delta >0$
with the property:
for any Cantor set $K\subset \mathbb R^N$
such that $d_H(f(\mathcal C) , K ) <\delta $
there exists an embedding
$g\in U (f,\varepsilon )$
with $g(\mathcal C) = K$.

In this proof, the notation of the form
$A=A_1\sqcup \ldots \sqcup A_n$
is used for a union that happens to be disjoint;
that is, $A_i\cap A_j = \emptyset $ for each
$i\neq j$, and 
$A=A_1\cup \ldots \cup A_n$.
Decompose
$f(\mathcal C) = X_1 \sqcup\ldots\sqcup X_s$,
where each $X_i$ is a Cantor set
with
$\diam X_i < \frac{\varepsilon }{3}$.
Take a positive number $\delta $
such that
$$\delta < \frac13 \min 
\{ d(X_i , X_j )
\mid
i,j=1,\ldots , s; \ 
i\neq j
  \}
\quad
\text{and}
\quad
\delta <\frac{\varepsilon }{3}.
$$
The neighborhoods
 $O(X_i , \delta )$, $i=1,\ldots , s$
are pairwise disjoint,
and 
\linebreak
$\diam O(X_i , \delta ) \leqslant \diam X_i + 2\delta
<\varepsilon $ for each $i$.
Let $K\subset \mathbb R^N$
be any set with
$d_H(f(\mathcal C) , K) < \delta $.
We have
$K\subset 
O(f(\mathcal C), \delta ) =
O(X_1, \delta )\sqcup
\ldots \sqcup O(X_s,\delta )$;
and it can be easily verified
that
$O(X_i,\delta ) \cap K \neq\emptyset $
for each $i=1,\ldots , s$.

Now take any Cantor set $K\subset \mathbb R^N$
 such that
$d_H(f(\mathcal C) , K) < \delta $.
For each $i=1,\ldots , s$
the intersection
$O(X_i,\delta ) \cap K  = : K_i$
is a Cantor set.
We have
$\mathcal C = f^{-1} (X_{1}) \sqcup
\ldots\sqcup  f^{-1} (X_{s})$,
and each $f^{-1}(X_i)$ is a Cantor set.
For any $i=1,\ldots , s$,
fix a homeomorphism
$h_i : f^{-1} ( X_{i} ) \cong K_i$.
Define 
an embedding
$g:\mathcal C\to\mathbb R^N$
by 
$g|_{f^{-1}(X_i)} = h_i$; 
we get
$$
\mathcal C = f^{-1} (X_{1}) \sqcup
\ldots\sqcup  f^{-1} (X_{s})
\stackrel{g=h_1\sqcup \ldots \sqcup h_s}{\to }
K_1\sqcup\ldots \sqcup K_s = K
\subset \mathbb R^N .
$$
We have
$g(\mathcal C) = K$.
The inclusions
$$f(f^{-1} (X_i )) \cup
g(f^{-1} (X_i ))  = X_i \cup K_i 
\subset O(X_i , \delta ),
\quad i=1,\ldots , s$$
together with the inequality
$\diam O(X_i , \delta ) <\varepsilon $
imply that
$\rho (g,f)<\varepsilon $.
Thus $g$ is the required embedding.

Proof of Statement \ref{easier}.
1)
Consider the map
$\Pi : 
\mathcal C (\mathbb R^N)\to
\mathcal K (L)$,
$X\mapsto p_{L} (X)$.
We have
$d_H(X_1,X_2)\geqslant d_H(\Pi (X_1),\Pi (X_2))$
for each $X_1,X_2\in \mathcal C(\mathbb R^N)$,
hence $\Pi $ is continuous
(see Lemma \ref{prelim}).

The map $\Pi $ is surjective.
To prove this, take a non-empty compactum $Y\subset  L $. There exists
a continuous surjection
$f: \mathcal C \to Y$ \cite[4.18]{Kechris};
let $\Gamma (f)\subset \mathcal C\times Y$
be its graph.
Decompose 
the 
orthogonal complement
$L^\bot $ of $L$ in
$\mathbb R^N$
as
$L^\bot = V^\bot \oplus V$, 
where $V$ is a one-dimensional subspace of $L^\bot $.
Place 
 a copy of 
 $\Gamma (f)$ 
 in $\mathbb R^N$
 as follows:
$$ 
\Gamma (f)\subset \mathcal C\times Y
\subset \mathbb R \times L
\stackrel{\xi \times \id _L}{\to  }
\{ 0 \} \oplus V \oplus L\subset
L^\bot \oplus L =
\mathbb R^N ,
$$
where $\xi : \mathbb R\to \{0\}\oplus V$ is any isomorphism.
The image of
$\Gamma (f)$ under this embedding
is a Cantor set in $\mathbb R^N$
whose projection to $L $
coincides with $Y$.

Our aim is to show that the set
$\mathcal P_1 = \Pi  ^{-1} (\mathcal C (L ))$
is a (dense) $G_\delta $ subset of
$\mathcal C (\mathbb R^N)$.
Recall that $\mathcal C(L )$
is a dense $G_\delta $ subset of
$\mathcal K(L)$,
hence
$\Pi  ^{-1} (\mathcal C (L ))$
is 
a $G_\delta $ subset of
$\mathcal C (\mathbb R^N)$.
It remains to prove
the density of
$\Pi ^{-1} (\mathcal C (L )) $
in 
$\mathcal C (\mathbb R^N)$.
Take any Cantor set $K\subset \mathbb R^N$
and any number $\varepsilon >0$.
We will construct a Cantor set $X \subset \mathbb R^N$
such that  $d_H(K,X) < \varepsilon  $
and $p_{L } X$ is a Cantor set.

 Take a  finite set
$A= \{ a_1,\ldots , a_t\}
\subset\mathbb R^N$ such that
$d_H(K,A)<\frac{\varepsilon }{2}$
(Lemma~\ref{finite-appr}).
Let $r$ be a positive number such that $r < \frac12 \min \{ d(x,y) \mid x,y\in A,\ x\neq y\}$
and $r<\frac{\varepsilon }{2}$.
For each $j=1,\ldots , t$
take a line segment
$\ell _j \subset O (a_j , r)$
parallel to $L $,
and a Cantor set $X_j \subset \ell _j$.
The union $X=X_1\cup \ldots \cup X_s$
is the desired Cantor set.
In fact, $p_L X$ is a Cantor set, and
$$d_H(K,X) \leqslant
d_H(K,A) + d_H(A,X) <
\frac{\varepsilon }{2} +\frac{\varepsilon }{2}  = \varepsilon  .$$

2) As in 1), 
the
map 
$\widetilde \Pi : \mathcal K (\mathbb R^N) 
\to
\mathcal K (L)$,
$X\mapsto p_L X$,
is continuous.
Hence the set
$\mathcal P_2 = \widetilde \Pi ^{-1} (\mathcal C(L))$
is a $G_\delta $ subset of $\mathcal K(\mathbb R^N)$;
it remains to verify its density.
By 1),
$\Pi ^{-1}(\mathcal C(L))$
is dense in 
$\mathcal C(\mathbb R^N)$.
Recall that
$\mathcal C(\mathbb R^N)$
is dense in 
$\mathcal K (\mathbb R^N)$
(see the references 
straight after Theorem~\ref{main}).
Clearly
$\Pi ^{-1}(\mathcal C(L))
\subset 
\widetilde \Pi ^{-1}(\mathcal C(L))$;
the result follows.

3)
We have 
$\mathcal P_3 =
\mathcal E^{-1} (\mathcal P_1)$,
where 
$\mathcal E :
\Emb (\mathcal C,\mathbb R^N)
\to 
\mathcal C(\mathbb R^N)$
takes $f\mapsto f(\mathcal C)$.
The proposition follows 
from 1) together with Lemma~\ref{openness}.

\end{document}